\documentclass[11pt,reqno]{amsart}
\usepackage{amsmath,amssymb,latexsym,textcomp,mathrsfs}
\usepackage[all]{xy}
\usepackage{graphicx}
\usepackage{bm,amsmath, amsthm, amssymb, amsfonts}
\usepackage{times}
\usepackage[english]{babel}

\setbox0=\hbox{$+$}
\newdimen\plusheight
\plusheight=\ht0
\def\+{\;\lower\plusheight\hbox{$+$}\;}  

\setbox0=\hbox{$-$}
\newdimen\minusheight
\minusheight=\ht0
\def\-{\;\lower\minusheight\hbox{$-$}\;}

\setbox0=\hbox{$\cdots$}
\newdimen\cdotsheight
\cdotsheight=\plusheight
\def\cds{\lower\cdotsheight\hbox{$\cdots$}}

\numberwithin{equation}{section}
\theoremstyle{plain}

  \newenvironment{nouppercase}{%
   \renewcommand{\uppercasenonmath}[1]{}}{}
	
	 \newcommand{\Keywords}[1]{\par\noindent
   {\small{\textbf{Keywords and phrases}}: #1}}
   
   \newcommand{\AMS}[1]{\par\noindent
   {\small{\textbf{AMS Subject Classification (2010)}}: #1}}

 \author{Amar Kumar Banerjee$^{*}$}
 \author{Apurba Banerjee}

\begin{document}
 
\title{$\textit{I}$-convergence classes of sequences and nets in topological spaces}
   
\begin{abstract}
In this paper we have used the idea of $I$-convergence of sequences and nets to study certain conditions of convergence in a topological space. It has been shown separately that a class of sequences and a class of nets in a non-empty set $X$ which are respectively called $I$-convergence class of sequences and $I$-convergence class of nets satisfying these conditions generate a topology on $X$. Further we have correlated the classes of $I$-convergent sequences and nets with respect to these topologies with the given classes which satisfy these conditions.
\end{abstract}

\begin{nouppercase}
\maketitle
\end{nouppercase}

\Keywords{ideal, filter, net, $\textit{I}$-convergence, $I$-convergence class, $\textit{I}$-cluster point, $I$-limit space.}  \\
\AMS{Primary 54A20, Secondary 40A35.}
\\

\let\thefootnote\relax\footnotetext{
$^{*}$Corresponding Author  \\
Amar Kumar Banerjee \\
Department of Mathematics, The University of Burdwan  \\
Golapbag, Burdwan-713104, West Bengal, India. \\
e-mail: akbanerjee@math.buruniv.ac.in, akbanerjee1971@gmail.com  \\
\\
Apurba Banerjee   \\
Department of Mathematics, The University of Burdwan  \\
Golapbag, Burdwan-713104, West Bengal, India. \\
e-mail: apurbabanerjeemath@gmail.com
}

\section{\bf Introduction and background}
The concept of convergence of a sequence of real numbers was extended to statistical convergence independently by H.Fast \cite{FH} and I.J.Schoenberg \cite{SIJ} as follows: 

If $K$ is a subset of the set of all natural numbers $\mathbb{N}$ then natural density of the set $K$ is defined by $d(K)=\lim_{n\rightarrow\infty}\frac {\left\vert K_{n}\right\vert}{n}$ if the limit exits (\cite{HH},\cite{NI}) where $\left\vert K_{n}\right\vert$ stands for the cardinality of the set $K_{n}=\left\{k \in K:k \leq n \right\}.$ 

A sequence $\left\{x_{n}\right\}$ of real numbers is said to be \textit{statistically convergent} to $\ell$ if for every $\varepsilon>0$ the set 
$$K(\varepsilon)=\left\{k \in \mathbb{N}:\left\vert x_{k}-\ell \right\vert \geq \varepsilon \right\}$$
has natural density zero (\cite{FH},\cite{SIJ}).

This idea of statistical convergence of real sequence was generalized to the idea of $\textit{I}$-convergence of real sequences (\cite{KP},\cite{KPMM}) using the notion of ideal $\textit{I}$ of subsets of the set of natural numbers. Several works on $\textit{I}$-convergence and on statistical convergence have been done in (\cite{BV},\cite{DK},\cite{KP},\cite{KPMM},\cite{LD1},\cite{MM}).

The idea of $\textit{I}$-convergence of real sequences coincides with the idea of ordinary convergence if $\textit{I}$ is the ideal of all finite subsets of $\mathbb{N}$ and with the statistical convergence if $\textit{I}$ is the ideal of subsets of $\mathbb{N}$ of natural density zero.
Later B.K. Lahiri and P. Das (\cite{LD2}) extended the idea of $\textit{I}$-convergence to an arbitrary topological space and observed that the basic properties are preserved also in a topological space. 
They also introduced (\cite{LD3}) the idea of $\textit{I}$-convergence of nets in a topological space and examined how far it affects the basic properties.      \\
The study of Moore-Smith convergence of sequences and nets (\cite{KJL}) deals with the construction of a topology on a given non-void set $X$ as follows:

Let $\mathcal{C}$ be a class consisting of pairs $(S,s)$ where $S$ is a net in $X$ and $s$ is a point of $X$. Then $\mathcal{C}$ is called a convergence class for $X$ if and only if it satisfies the conditions \textbf{(a)} to \textbf{(d)} given below. For convenience, we say that $S$ converges $(\mathcal{C})$ to $s$ or lim$_{n}S_{n}$=$s$ $(\mathcal{C})$ if and only if $(S,s)\in \mathcal{C}$.   \\
\\
\textbf{(a)}  If $S$ is a net such that $S_{n}$=$s$ for each $n$, then $S$ converges $(\mathcal{C})$ to $s$.  \\
\textbf{(b)}  If $S$ converges $(\mathcal{C})$ to $s$, then so does each subnet of $S$.   \\
\textbf{(c)}  If $S$ does not converge $(\mathcal{C})$ to $s$, then there is a subnet of $S$, no subnet of which converges $(\mathcal{C})$ to $s$.    \\
\textbf{(d)}  (\textbf{Theorem on iterated limits}) Let $D$ be a directed set and let $E_{m}$ be a directed set for each $m$ in $D$. Let $F$ be the product $D\times (\times\{E_{m}: m\in D\})$ and for $(m,f)$ in $F$ let $R(m,f)$=$(m,f(m))$. If $lim_{m}lim_{n} S(m,n)$=$s$ $(\mathcal{C})$, then $S\circ R$ converges $(\mathcal{C})$ to $s$. 

Indeed if $S$ is a net in a topological space $(X,\tau)$, then convergence of $S$ with respect to the topology $\tau$ implies all the conditions listed above and in turn a convergence class $\mathcal{C}$ determines a topology $\sigma$ on $X$ such that $(S,s)\in \mathcal{C}$ if and only if $S$ converges to $s$ relative to this topology $\sigma$.
The study of convergence class of sequences in $X$ and construction of topology is almost similar to that of convergence class of nets.   \\

Here we have used the idea of $I$-convergence of sequences and nets to study certain conditions of convergence of sequences and nets which are in turn sufficient to determine a topology on a given non-void set $X$. Also we have obtained a correlation between the given classes of $I$-convergent sequences and nets satisfying these conditions and the classes of $I$-convergent sequences and nets with respect to the topologies generated by the given classes of $I$-convergent sequences and nets.     \\
\section{\bf $I$-convergence class of sequences and $I$-limit space}
First we recall the following definitions. \\

\textbf{Definition 2.1}(\cite{KK}) If $\textit{X}$ is a non-void set then a family of sets $\textit{I}\subset 2^{\textit{X}}$ is called an \textit{ideal} if \\
(i) $A,B \in \textit{I}$ implies $A \cup B \in \textit{I}$ and \\
(ii) $A\in \textit{I},B\subset A$ imply $B\in \textit{I}.$ \\

The ideal $I$ is called \textit{nontrivial} if $\textit{I} \neq \left\{\emptyset\right\}$ and $\textit{X}\notin \textit{I}$. \\

\textbf{Definition 2.2}(\cite{KK}) A non-empty family \textit{F} of subsets of a non-void set \textit{X} is called a \textit{filter} if \\
(i)  $\emptyset\notin \textit{F}$ \\
(ii) $A,B\in \textit{F}$ implies $A\cap B\in \textit{F}$ and \\
(iii) $A\in\textit{F},A\subset B$ imply $B\in \textit{F}.$ \\

If \textit{I} is a nontrivial ideal on \textit{X} then $\textit{F}=\textit{F}(\textit{I})=\left\{A\subset \textit{X}:X-A \in \textit{I}\right\}$ is clearly a filter on \textit{X} and conversely. \\

A nontrivial ideal \textit{I} is called \textit{admissible} if it contains all the singleton sets. Several examples of nontrivial admissible ideals may be seen in \cite{KP}. \\

Let $(\textit{X},\tau)$ be a topological space and \textit{I} be a nontrivial ideal of $\mathbb{N}$, the set of all natural numbers. \\

\textbf{Definition 2.3}(\cite{LD2}) A sequence $\left\{x_{n}\right\}$ in \textit{X} is said to be $\textit{I}$-convergent to $x_{0}\in \textit{X}$ if for any non-empty open set $U$ containing $x_{0}$, $\left\{n\in \mathbb{N}:x_{n}\notin U\right\}\in \textit{I}.$ \\

In this case, $x_{0}$ is called an $\textit{I}$-limit of $\left\{x_{n}\right\}$ and written as $x_{0}$=$\textit{I}$-lim $x_{n}.$ \\

\textbf{Remark 2.1} If \textit{I} is an admissible ideal then ordinary convergence implies $\textit{I}$-convergence and if \textit{I} does not contain any infinite set then converse is also true. \\

\textbf{Definition 2.4}(\cite{LD2}) An element $y\in X$ is said to be an $I$-cluster point of a sequence $\{x_{n}\}$ of elements of $X$ if for every non-empty open set $U$ containing $y$ we have the set $\{n\in \mathbb{N}: x_{n}\in U\} \notin I$.      \\

We prove below some properties of a convergent sequence in a topological space which remains invariant in case of $I$-convergence of a sequence in a topological space. 

Throughout $(X,\tau)$ stands for a topological space and $I$ a nontrivial admissible ideal of $\mathbb{N}$.    \\

\textbf{Theorem 2.1} Let $(X,\tau)$ be a topological space. Then the following conditions hold   \\
\textbf{C(1):} For any point $x_{0}\in X$ the sequence $\{x_{0},x_{0},x_{0},\ldots\}$ is $I$-convergent to $x_{0}$.   \\
\textbf{C(2):} Addition of a finite number of terms to a sequence affects neither its $I$-convergence nor its $I$-limit.    \\
\textbf{C(3):} If  a sequence $\{x_{n}\}$ in $(X,\tau)$ is $I$-convergent to $x_{0}\in X$ then every subsequence of it is $I$-convergent to the same $I$-limit $x_{0}$.   

\begin{proof}
Since for any non-empty open set $U$ containing $x_{0}$ we have $\{n\in \mathbb{N}: x_{n}\notin U\}$=$\emptyset \in I$, the property C(1) holds.   \\ 
Let $\{x_{n}\}$ be a sequence in $(X,\tau)$ which is $I$-convergent to $x_{0}\in X$. Now let finite number of points say $y_{1},y_{2},\ldots,y_{r}$ be included into the sequence $\{x_{n}\}$ and let us denote the new sequence by $\{z_{n}\}$. Then for any non-empty open set $U$ containing $x_{0}$ we have $\{n\in \mathbb{N}: z_{n}\notin U\}$=$\{n\in \mathbb{N}: x_{n}\notin U\} \cup \{n\in \mathbb{N}: y_{n}\notin U\}$. Now the first set on the right hand side belongs to $I$ and the second set being a finite set also belongs to $I$, since $I$ is an admissible ideal. Thus C(2) holds.    \\
The proof of C(3) follows from the definitions of an ideal and $I$-convergence of a sequence.
\end{proof}

The following two properties hold in a topological space in case of ideal convergence.    \\

\textbf{Theorem 2.2} Let $G$ be an open set in $(X,\tau)$. Then no sequence lying in $X-G$ has any $I$-limit in $G$.

\begin{proof}
If possible let $\{x_{n}\}$ be a sequence in $X-G$ which is $I$-convergent to $x_{0}\in G$. Since $G$ is an open set containing $x_{0}$, we must have by definition of $I$-convergence that the set $\{n\in \mathbb{N}: x_{n}\notin G\}\in I$ i.e., $\mathbb{N}\in I$, which leads to a contradiction, since $I$ is a non-trivial ideal. Hence the proof follows.
\end{proof}

\textbf{Theorem 2.3} If $F$ be a closed set in $(X,\tau)$ then every $I$-convergent sequence lying in $F$ has all its $I$-limits in $F$.

\begin{proof}
The proof is similar to the proof of Theorem 2.2. 
\end{proof}

We shall now show that a topology can be generated in terms of a class of sequences satisfying the conditions C(1),C(2),C(3) of Theorem 2.1. In fact the open sets are determined by the conditions above.   \\

\textbf{Theorem 2.4} Let $X$ be a given non-void set and let a class of infinite sequences $\Omega$ over X, be distinguished whose members are called '$I$-convergent sequences', and let each $I$-convergent sequence be associated with an element of $X$ which is called an '$I$-limit' of the sequence subject to the conditions C(1),C(2),C(3) of Theorem 2.1.

Let now, a subset $G$ of $X$ be called an open set, if and only if no sequence lying in $X-G$ has any $I$-limit in $G$. Then the collection of all such open sets thus obtained forms a topology $\sigma$ on $X$.  

\begin{proof}
Clearly $\emptyset$ and $X$ are open sets.    \\
Let $\{G_{\lambda}\}_{\lambda\in \Lambda}$ be a collection of open sets and $G$=$\bigcup_{\lambda\in \Lambda} G_{\lambda}$, where $\Lambda$ is an indexing set. If possible let $G$ be not an open set. Then there exists an $I$-convergent sequence $\{x_{n}\}$ in $X-G$ which has an $I$-limit say $x_{0}$ in $G$. Then $x_{0}\in G_{\lambda}$ for some $\lambda \in \Lambda$. So the sequence $\{x_{n}\}$ lying in $X-G_{\lambda}$ has an $I$-limit $x_{0}\in G_{\lambda}$ which is impossible, since $G_{\lambda}$ is an open set. Hence $G$ must be an open set.      \\
Let $G,H$ be two open sets. If possible let $G\cap H$ be not an open set. Then there exists a sequence $\{x_{n}\}\subset X-(G\cap H)$ with an $I$-limit say $x_{0}$ in $G\cap H$. Now since $x_{0}\in G\cup H$ and $G\cup H$ is open, only a finite number of terms can lie outside $G\cup H$. Again since $G\cup H$=$(G-H)\cup (G\cap H)\cup (H-G)$ and since $\{x_{n}\}$ lies wholly outside $G\cap H$ we must have either $G-H$ or $H-G$ contains infinite number of terms of $\{x_{n}\}$. Let us suppose that $G-H$ contains infinite number of terms of $\{x_{n}\}$. Then we get a subsequence of $\{x_{n}\}$ lying wholly outside $H$ which is $I$-convergent to a point of $H$. But this contradicts that $H$ is an open set. So $G\cap H$ must be an open set.
\end{proof} 

The topology $\sigma$ defined as above is called $I$-convergence topology on $X$ and $(X,\sigma)$ is called $I$-limit space.   \\

It should however be noted that the family of all $I$-convergent sequences in $I$-limit space $(X,\sigma)$ need not be identical with the given family $\Omega$ stated in Theorem 2.4. However the following results are true      \\

\textbf{Theorem 2.5} Let $(X,\sigma)$ be an $I$-limit space with $\Omega$ as the given collection of $I$-convergent sequences on $X$. If $\Gamma$ is the collection of all $I$-convergent sequences with respect to the topology $\sigma$ on $X$, then $\Omega\subset \Gamma$.

\begin{proof}
Let $\{x_{n}\}$ be a sequence in $\Omega$ with '$I$-limit' $x_{0}$. Let $G$ be an open set in $(X,\sigma)$ containing $x_{0}$. Then only a finite number of terms can possibly lie outside $G$, because otherwise an infinite subsequence of $\{x_{n}\}$ lying outside $G$ would have an $I$-limit $x_{0}$ in $G$ which leads to a contradiction, since $G$ is an open set. So $\{n\in \mathbb{N}: x_{n}\notin G\}\in I$, since $I$ is a nontrivial admissible ideal. Therefore $\{x_{n}\}$ is $I$-convergent to $x_{0}$ with respect to the topology $\sigma$. Hence $\Omega \subset \Gamma$. 
\end{proof}

\textbf{Theorem 2.6} Let $\Gamma$ be the family of $I$-convergent sequences in a topological space $(X,\tau)$ and let $\tau^{\prime}$ be the $I$-convergence topology on $X$ determined by the family $\Gamma$. Then $\tau\subset \tau^{\prime}$. 

\begin{proof}
Let $G\in \tau$ and $\{x_{n}\}$ be a sequence in $\Gamma$ which is $I$-convergent to $x_{0}\in G$. Then the set $\{n\in \mathbb{N}: x_{n}\notin G\}\in I$. So $\{x_{n}\}$ cannot lie wholly in $X-G$, since $I$ is a nontrivial ideal. So we conclude that no sequence lying wholly in $X-G$ can be $I$-convergent to a point in $G$. Hence $G$ becomes $\tau^{\prime}$-open. Thus $\tau \subset \tau^{\prime}$.
\end{proof}
\hspace{0.5in}
\section{\bf $I$-convergence classes of nets and $I$-convergence topology}  
The following definitions are widely known.  \\

\textbf{Definition 3.1} A binary relation $\geq$ directs a set $D$ if $D$ is non-void and  \\
\textbf{(a)} $m\geq m$ for each $m\in D$;   \\
\textbf{(b)} if $m$,$n$ and $p$ are members of $D$ such that $m\geq n$ and $n\geq p$, then $m\geq p$; and     \\
\textbf{(c)} if $m$ and $n$ are members of $D$, then there is a member $p$ of $D$ such that $p\geq m$ and $p\geq n$. 
\\
The pair $(D,\geq)$ is called a \textit{directed set}.   \\

\textbf{Definition 3.2} Let $(D,\geq)$ be a directed set and let $X$ be a non-void set. A mapping $S:D\rightarrow X$ is called a \textit{net} in $X$ denoted by $\{S_{n}:n\in D\}$ or simply by $\{S_{n}\}$ when the set $D$ is clear.     \\

\textbf{Definition 3.3} Let $(D,\geq)$ be a directed set and $\{S_{n}:n\in D\}$ be a net in $X$. A net $\{T_{m}:m\in E\}$ where $(E,\succ)$ is a directed set is said to be a \textit{subnet} of $\{S_{n}:n\in D\}$ if and only if there is a function $\theta$ on $E$ with values in $D$ such that   \\
\textbf{(a)} $T$=$S\circ \theta$, or equivalently, $T_{i}$=$S_{\theta_{i}}$ for each $i\in E$; and  \\
\textbf{(b)} for each $m$ in $D$ there is $n$ in $E$ with the property that, if $p\in E$ and $p\geq n$, then $\theta_{p}\geq m$.
\\

Throughout our discussion $(X,\tau)$ will denote a topological space and $I$ will denote a nontrivial ideal of a directed set $D$ unless otherwise stated. \\

For $n\in D$ let $M_{n}$=$\{k\in D: k\geq n\}$. Then the collection $F_{0}$=$\{A\subset D: A\supset M_{n}$ for some $n\in D \}$ forms a filter in $D$. Let $I_{0}$=$\{B\subset D: D-B\in F_{0}\}$. Then $I_{0}$ is also a nontrivial ideal of $D$.   \\

\textbf{Definition 3.4} (\cite{LD3}) A nontrivial ideal $I$ of $D$ will be called \textit{$D$-admissible} if $M_{n}\in F(I)$ for all $n\in D$, where $F(I)$ is the filter associated with the ideal $I$ of $D$.      \\

\textbf{Definition 3.5} (\cite{LD3}) A net $\{S_{n}:n\in D\}$ in $X$ is said to be $I$-convergent to $x_{0}\in X$ if for any open set $U$ containing $x_{0}$ the set $\{n\in D: S_{n}\notin U\}\in I$.   \\

Symbolically we write $I$-lim $S_{n}$=$x_{0}$ and we say that $x_{0}$ is an $I$-limit of the net $\{S_{n}\}$.  \\

\textbf{Remark 3.1} If $I$ is $D$-admissible, then convergence of net in the topology $\tau$ implies $I$-convergence and the converse holds if $I$=$I_{0}$. In other words, $I_{0}$-convergence implies net convergence.  \\

\textbf{Definition 3.6} (\cite{LD3}) A point $y\in X$ is called an $I$-cluster point of a net $\{S_{n}:n\in D\}$ if for every open set $U$ containing $y$, $\{n\in D: S_{n}\in U\}\notin I$.    \\ 

The following result is very useful.         \\

\textbf{Theorem 3.1} (\cite{LD3}) Let $\{S_{n}:n\in D\}$ be a net in a topological space $(X,\tau)$ and $I$ be a nontrivial ideal of $D$. Then $x_{0}\in X$ is an $I$-cluster point of $\{S_{n}\}$ if and only if $x_{0}\in \overline{A_{T}}$ for every $T\in F(I)$, where $A_{T}$=$\{x\in X: x$=$S_{t}$ for $t\in T\}$ and $F(I)$ is the filter associated with the ideal $I$ of $D$. Here bar denotes the closure in $(X,\tau)$.  
\\

In case of ideal convergence of a subnet of a net in a topological space $(X,\tau)$ the following results hold.    \\

\textbf{Theorem 3.2} Let $\{S_{n}:n\in D\}$ be a net in a topological space $(X,\tau)$ and $I_{D}$ be a nontrivial ideal of $D$. Let $\{T_{m}: m\in (E,\succ)\}$ be a subnet of $\{S_{n}:n\in D\}$ and $I_{E}$=$\{A\subset E: \theta(A)\in I_{D}\}$ where $\theta: E\rightarrow D$ is the function associated with $\{T_{m}: m\in (E,\succ)\}$ to be a subnet of $\{S_{n}:n\in D\}$. Then $I_{E}$ is an ideal of $E$ and if $\{S_{n}: n\in D\}$ be $I_{D}$-convergent to $x_{0}\in X$ and $I_{E}$ is a nontrivial ideal of $E$ then $\{T_{m}: m\in E\}$ is $I_{E}$-convergent to $x_{0}$. 

\begin{proof}
Since $\{T_{m}: m\in (E,\succ)\}$ is a subnet of $\{S_{n}:n\in D\}$ so $\theta: E\rightarrow D$ is a function such that $T_{m}$=$S\circ \theta(m)$ for all $m\in E$ i.e., $T_{m}$=$S_{\theta_{m}}$ for all $m\in E$. Now since $\{S_{n}: n\in D\}$ is $I_{D}$-convergent to $x_{0}$ so for every open set $U$ containing $x_{0}$, the set $\{n\in D: S_{n}\notin U\}\in I_{D}$. If possible let $\{T_{m}: m\in E\}$ be not $I_{E}$ convergent to $x_{0}$. Then there exists an open set $V$ containing $x_{0}$ such that the set $M$=$\{m\in E: T_{m}$=$S_{\theta_{m}}\notin V\}\notin I_{E}$. Then by definition of the ideal $I_{E}$, the set $\theta(M)$=$\{\theta(m)\in D: S_{\theta_{m}}\notin V\}\notin I_{D}$. But since $\theta(M)$=$\{\theta(m)\in D: S_{\theta_{m}}\notin V\} \subset \{n\in D: S_{n}\notin V\}\in I_{D}$, by definition of ideal we get that $\theta(M)\in I_{D}$. Thus we arrive at a contradiction. Hence we must have the set $M$=$\{m\in E: T_{m}$=$S_{\theta_{m}}\notin V\}\in I_{E}$. Therefore the result follows.
\end{proof}

\textbf{Theorem 3.3} Let $\{S_{n}:n\in D\}$ be a net in a topological space $(X,\tau)$. Let $I_{D}$ be a $D$-admissible ideal of $D$. Let $\{S_{n}\}$ be not $I_{D}$-convergent to a point $x_{0}\in X$. Then there exists a subnet of $\{S_{n}\}$, no subnet of which is ideal convergent to $x_{0}$ with respect to any nontrivial ideal.

\begin{proof}
Let the net $\{S_{n}:n\in D\}$ be not $I_{D}$-convergent to $x_{0}\in X$. Then there exists a non-empty open set $U$ containing $x_{0}$ such that the set $A$=$\{n\in D: S_{n}\notin U\}\notin I_{D}$ and consequently the set $A^{c}$=$\{p\in D: S_{p}\in U\}\notin F(I_{D})$, where $F(I_{D})$ is the filter on $D$ associated with the ideal $I_{D}$. Since $I_{D}$ is $D$-admissible so $F(I_{D})$ contains all sets of the form $M_{n}$=$\{m\in D: m\geq n\}$ for each $n\in D$. Again since $A^{c}\notin F(I_{D})$, we conclude that for no $n\in D$, $M_{n}\subset A^{c}$. Consequently for each $n\in D$ there is some $m\in M_{n}$ such that $m\notin A^{c}$ and so $S_{m}\notin U$. Let $B_{n}$=$\{m\in M_{n}: m\notin A^{c}\}$ and $M$=$\bigcup_{n\in D} B_{n}$. Then clearly $M$ is a cofinal subset of $D$. Let $\{T_{r}:r\in M\}$ be a subnet of $\{S_{n}:n\in D\}$. Then we see that $T_{r}\notin U$, for all $r\in M$. Now if $\{K_{p}:p\in E\}$ be a subnet of $\{T_{r}:r\in M\}$ where $(E,\succ)$ is a directed set and $I_{E}$ be any nontrivial ideal of $E$ then we note that $\{K_{p}:p\in E\}$ is not $I_{E}$-convergent to $x_{0}$, since for the open set $U$ containing $x_{0}$ we have $\{p\in E: K_{p}\notin U\}$=$E\notin I_{E}$. Hence the result follows.
\end{proof} 

We now recall the definition of product directed set.   \\

Suppose that for each member $a$ of a set $A$ we are given a directed set $(D_{a},>_{a})$, where $A$ is an indexing set. The cartesian product $\times\{D_{a}:a\in A\}$ is the set of all functions $d$ on $A$ such that $d_{a}$(=$d(a)$) is a member of $D_{a}$ for each $a$ in $A$. The product directed set is $(\times\{D_{a}:a\in A\},\geq)$ where, if $d$ and $e$ are members of the product $\times\{D_{a}:a\in A\}$ then $d\geq e$ if and only if $d_{a}>_{a}e_{a}$ for each $a$ in $A$. The product order is $\geq$.      \\

\textbf{Definition 3.7} Let $D$ be a directed set and for each $m\in D$, let $E_{m}$ be a directed set. Consider a function $S$ to a topological space $(X,\tau)$ such that $S(m,n)$ is defined whenever, $m\in D$, $n\in E_{m}$. Let $I_{D}$ be a nontrivial ideal of $D$ and $I_{E_{m}}$ be a nontrivial ideal of $E_{m}$ for each $m\in D$. We say that $I_{D}$-lim$_{m}I_{E_{m}}$-lim$_{n}S(m,n)$=$x_{0}\in X$ if for any non-empty open set $U$ containing $x_{0}$ we have the set $\{m\in D: I_{E_{m}}$-lim$_{n}S(m,n)\notin U\}\in I_{D}$.    \\

\textbf{Theorem 3.4} (\textbf{THEOREM ON ITERATED $I$-LIMIT}) Let $D$ be a directed set, let $E_{m}$ be a directed set for each $m$ in $D$. Let $\mathcal{F}$ be the product $D\times (\times\{E_{m}:m\in D\})$ and for $(m,f)$ in $\mathcal{F}$ let $R(m,f)$=$(m,f(m))$. Let $I_{D}$ be a nontrivial ideal of $D$ and for each $m\in D$ and let $I_{E_{m}}$ be a nontrivial ideal of $E_{m}$. Let $I_{\mathcal{F}}$ be a nontrivial ideal of $\mathcal{F}$ defined as follows:  \\
A subset $H\subset \mathcal{F}$ belongs to $I_{\mathcal{F}}$ if and only if $H$=$H_{1}\cup H_{2}$, where $H_{2}$ may be empty set and $H_{1}$, $H_{2}$ are such that   
$$H_{1}^{\prime}=\{m\in D: (m,f)\in H_{1}\}\in I_{D},$$     
$$H_{2}^{\prime}=\{p\in D: (p,g)\in H_{2}\}\notin I_{D}$$    
and the set    
$$\{g(p): (p,g)\in H_{2}\}\in I_{E_{p}}$$ for each fixed $p\in H_{2}^{\prime}$ if $H_{2}\neq \emptyset$. If $S(m,n)$ is a member of a topological space $(X,\tau)$ for each $m$ in $D$ and $n$ in $E_{m}$ then $S\circ R$ is $I_{\mathcal{F}}$-convergent to $I_{D}$-lim$_{m}I_{E_{m}}$-lim$_{n}S(m,n)$ whenever this iterated limit exists. 

\begin{proof}
Let $I_{D}$-lim$_{m}I_{E_{m}}$-lim$_{n}S(m,n)$=$x_{0}$ and $U$ be an open set containing $x_{0}$. Then the set $\{m\in D: I_{E_{m}}$-lim$_{n}S(m,n)\notin U\}\in I_{D}$. Consequently, the set $A$=$\{m\in D: I_{E_{m}}$-lim$_{n}S(m,n)\in U\}\notin I_{D}$. So for each $m\in A$, the set $A_{m}$=$\{n\in E_{m}: S(m,n)\notin U\}\in I_{E_{m}}$. Consequently, for each $m\in A$, the set (say) $B_{m}$=$\{t\in E_{m}: S(m,t)\in U\}\notin I_{E_{m}}$. For each $m\in A$, let $C_{m}$=$\{f\in (\times\{E_{m}:m\in D\}): f(m)\in B_{m}\}$. Now we write $C$=$\{(m,f)\in \mathcal{F}: m\notin A,f\notin C_{m}\}\cup \{(m,f)\in \mathcal{F}: m\notin A,f\in C_{m}\}\cup \{(m,f)\in \mathcal{F}: m\in A,f\notin C_{m}\}$. Let $P$=$\{(m,f)\in \mathcal{F}: m\notin A,f\notin C_{m}\}$, $Q$=$\{(m,f)\in \mathcal{F}: m\notin A,f\in C_{m}\}$ and $R$=$\{(m,f)\in \mathcal{F}: m\in A,f\notin C_{m}\}$. Then $C$=$P\cup Q\cup R$. Let $M$=$\{(m,f)\in \mathcal{F}: S\circ R(m,f)\notin U\}$=$\{(m,f)\in \mathcal{F}: S(m,f(m))\notin U\}$. Let us take a member $(p,g)\in \mathcal{F}$ such that $(p,g)\notin C$. Then $(p,g)\in \mathcal{F}$ such that $p\in A$ and $g\in C_{m}$. This implies $g(p)\in B_{m}$ which in turn implies that $S(p,g(p))\in U$ i.e, $S\circ R(p,g)\in U$ i.e., $(p,g)\notin M$. Hence $M\subset C$. Now we see that the sets $P^{\prime}$=$\{m\in D: (m,f)\in P\}\in I_{D}$ and $Q^{\prime}$=$\{m\in D: (m,f)\in Q\}\in I_{D}$, since $P^{\prime},Q^{\prime}\subset D-A\in I_{D}$. Note that the set $R^{\prime}$=$\{m\in D: (m,f)\in R\}$ may or may not belong to $I_{D}$. Now we can write the set $C$ as below:     \\
$C$=$C_{1}\cup C_{2}$, where $C_{1}$=$P\cup Q$ and $C_{2}$=$R$ if $R^{\prime}\notin I_{D}$, otherwise $C_{1}$=$P\cup Q\cup R$ and $C_{2}$=$\emptyset$. The case that $C_{2}$=$\emptyset$ is trivial. If $C_{2}$ is non-empty then the set $\{f(m): (m,f)\in R\}\notin B_{m}$ which implies that $\{f(m)\in E_{m}: (m,f)\in R\}\subset A_{m}$, since $S(m,f(m))\notin U$ in this case. But this implies that the set $\{f(m)\in E_{m}: (m,f)\in R\}\in I_{E_{m}}$ for each fixed $m \in R^{\prime}$=$\{m\in D: (m,f)\in C_{2}\}$. Hence the result follows.
\end{proof}

Let $X$ be a fixed non-empty set and $\mathcal{M}$ be the class consisting of all pairs $(S,x_{0})$ where $\{S_{n}:n\in D\}$ is a net in $X$ and $x_{0}$ is a point of $X$. Throughout our discussion we will consider the following facts:    \\
If $\{S_{n}:n\in (D,\geq)\}$ be a net in $X$ then $I_{D}$ will denote a nontrivial $D$-admissible ideal of $D$ and if $\{T_{m}: m\in (E,\succ)\}$ be a subnet of $\{S_{n}:n\in D\}$ then $I_{E}$ will denote a nontrivial ideal of $E$ defined by $I_{E}$=$\{A\subset E: \theta(A)\in I_{D}\}$ where $\theta: E\rightarrow D$ is a function as in Definition 3.3. Also $F(I_{D})$ will denote the filter on $D$ associated with the ideal $I_{D}$ of $D$. \\
We shall say that $\mathcal{M}$ is an ideal convergence class for $X$ if and only if it satisfies the following conditions \textbf{(a)} to \textbf{(d)} listed below. For convenience we say that $S$ is $I_{D}$-convergent$(\mathcal{M})$ to $x_{0}$ or that $I_{D}$-lim$_{m}S_{m}$=$x_{0}$$(\mathcal{M})$ if and only if $(S,x_{0})\in \mathcal{M}$.       \\
\\
\textbf{(a)}  If $\{S_{n}:n\in D\}$ be a net such that $S_{n}$=$x_{0}$ for all $n\in D$, then $\{S_{n}\}$ is $I_{D}$-convergent$(\mathcal{M})$ to $x_{0}$.      \\
\textbf{(b)}  If a net $\{S_{n}:n\in D\}$ is $I_{D}$-convergent$(\mathcal{M})$ to $x_{0}$, then every subnet $\{T_{m}:m\in E\}$ is $I_{E}$-convergent$(\mathcal{M})$ to $x_{0}$.       \\
\textbf{(c)}  If $\{S_{n}:n\in D\}$ is not $I_{D}$-convergent$(\mathcal{M})$ to $x_{0}$, then there is a subnet of $\{S_{n}\}$, no subnet of which is ideal convergent$(\mathcal{M})$ to $x_{0}$ with respect to any nontrivial ideal.     \\
\textbf{(d)}  (\textbf{THEOREM ON ITERATED $I$-LIMIT}) Let $D$ be a directed set, let $E_{m}$ be a directed set for each $m\in D$. Let $\mathcal{F}$ be the product $D\times (\times\{E_{m}:m\in D\})$ and for $(m,f)$ in $\mathcal{F}$, let $R$ be the net defined by $R(m,f)$=$(m,f(m))$. Let $I_{D}$ be a nontrivial ideal of $D$ and for each $m\in D$ let $I_{E_{m}}$ be a nontrivial ideal of $E_{m}$ and $I_{\mathcal{F}}$ be a nontrivial ideal of $\mathcal{F}$ as defined in Theorem 3.4. Let $S(m,n)$ be a member of $X$ whenever $m\in D$ and $n\in E_{m}$. Now if $I_{D}$-lim$_{m}I_{E_{m}}$-lim$_{n}S(m,n)$=$x_{0}$$(\mathcal{M})$ then $S\circ R$ is $I_{\mathcal{F}}$-convergent$(\mathcal{M})$ to $x_{0}$.    \\

Already we have shown that if a net $\{S_{n}: n\in (D,\geq)\}$ is $I_{D}$-convergent to a point $s$ in a topological space $(X,\tau)$ then the conditions \textbf{(a), (b), (c)} and \textbf{(d)} are satisfied. We now show that every ideal convergence class $\mathcal{M}$ determines a topology on $X$ for which $\{S_{n}\}$ is $I_{D}$-convergent with respect to this topology if $\{S_{n}\}$ is $I_{D}$-convergent$(\mathcal{M})$. The converse part is also true if an additional condition $(J)$ holds.     \\

\textbf{Theorem 3.5} Let $\mathcal{M}$ be an ideal convergence class for a non-empty set $X$, and for each subset $A$ of $X$ let $A^{cl}$ be the set of all points $x_{0}$ such that, for some net $\{S_{n}:n\in D\}$ in $A$, $\{S_{n}\}$ is $I_{D}$-convergent$(\mathcal{M})$ to $x_{0}$. Then '$cl$' is a closure operator and if $(S,x_{0})\in \mathcal{M}$ then $S$ is $I_{D}$-convergent to $x_{0}$ with respect to the topology associated with the closure operator.                                       \\
Conversely, $(S,x_{0})\in \mathcal{M}$ if $\{S_{n}:n\in D\}$ is $I_{D}$-convergent to $x_{0}$ with respect to the topology associated with the closure operator and if the following additional condition \textbf{(J)} holds:    \\
\textbf{(J)}: Let $\{S_{n}:n\in D\}$ be a net in $X$ and $\{T_{m}: m\in (E,\succ)\}$ be a subnet of $\{S_{n}:n\in D\}$. If $I_{D}$ be a $D$-admissible ideal of $D$ then $I_{E}$ is an $E$-admissible ideal of $E$.

\begin{proof}
We first prove that '$cl$' is a closure operator. Since a net is a function on a directed set and the set is non-void by definition, so $\emptyset^{cl}$ is void. In view of condition $\textbf{(a)}$ for each member $y_{0}$ of a set $A$ there is a net $\{S_{n}:n\in D\}$ defined by $S_{n}$=$y_{0}$ for all $n\in D$, which is $I_{D}$-convergent$(\mathcal{M})$ to $y_{0}$ and hence $A\subset A^{cl}$. If $x_{0}\in A^{cl}$ then by the definition of the operator '$cl$' we have $x_{0}\in (A\cup B)^{cl}$ and consequently $A^{cl}\subset (A\cup B)^{cl}$ for each set $B$. Therefore $A^{cl}\cup B^{cl}\subset (A\cup B)^{cl}$. To show the reverse inclusion, suppose that $\{S_{n}:n\in D\}$ is a net in $A\cup B$ and let $\{S_{n}:n\in D\}$ be $I_{D}$-convergent$(\mathcal{M})$ to $x_{0}$. If $D_{A}$=$\{n\in D: S_{n}\in A\}$ and $D_{B}$=$\{n\in D: S_{n}\in B\}$ then $D_{A}\cup D_{B}$=$D$. Hence either $D_{A}$ or $D_{B}$ is cofinal in $D$ and so either $\{S_{n}: n\in D_{A}\}$ or $\{S_{n}: n\in D_{B}\}$ is a subnet of $\{S_{n}: n\in D\}$ which is also $I_{D_{A}}$-convergent$(\mathcal{M})$ or $I_{D_{B}}$-convergent$(\mathcal{M})$ respectively to $x_{0}$, by virtue of the condition $\textbf{(b)}$. Hence we get that $x_{0}\in A^{cl}\cup B^{cl}$ and thus we have shown that $A^{cl}\cup B^{cl}$=$(A\cup B)^{cl}$. We now show that $(A^{cl})^{cl}$=$A^{cl}$. If $\{T_{m}:m\in D\}$ is a net in $A^{cl}$ 	which is $I_{D}$-convergent$(\mathcal{M})$ to '$t$', then for each $m\in D$, there is a directed set $E_{m}$ and a net $\{S(m,n): n\in E_{m}\}$ which is $I_{E_{m}}$-convergent$(\mathcal{M})$ to $T_{m}$. Now condition $\textbf{(d)}$ shows that there is a net $\{R_{(m,f)}: (m,f)\in D\times (\times \{E_{m}:m\in D\})\}$ which is $I_{\mathcal{F}}$-convergent$(\mathcal{M})$ to $t$ and consequently $t\in A^{cl}$, where $\mathcal{F}$=$D\times (\times \{E_{m}:m\in D\})$. Hence $(A^{cl})^{cl}$=$A^{cl}$.        \\

We now prove that ideal convergence $(\mathcal{M})$ is identical with the ideal convergence relative to the topology $\tau$ associated with the operator '$cl$'.           \\
First suppose that $\{S_{n}:n\in D\}$ is $I_{D}$-convergent$(\mathcal{M})$ to $x_{0}$ and $S$ is not $I_{D}$-convergent to $x_{0}$ relative to the topology $\tau$. Then there is an open set $U$ containing $x_{0}$ such that the set $M$=$\{n\in D: S_{n}\notin U\}\notin I_{D}$. So the set $K$=$D-M$=$\{n\in D: S_{n}\in U\}\notin F(I_{D})$. Now $I_{D}$ being $D$-admissible ideal of $D$, we have for each $r\in D$ the set $B_{r}$=$\{p\in D: p\geq r\}\in F(I_{D})$. Since $K\notin F(I_{D})$, $B_{r}$ is not a subset of $K$ for all $r\in D$. Hence for every $r\in D$, we can find some $p\in B_{r}$ such that $p\notin K$. Let us denote for each $r\in D$, the set $N_{r}$=$\{p\in B_{r}: p\notin K\}$ and $E$=$\bigcup_{r\in D}N_{r}$. Clearly $E$ is a cofinal subset of $D$. So $\{S_{n}:n\in E\}$ is a subnet of $\{S_{m}:m\in D\}$ and $S_{n}\notin U$ for all $n\in E$. Again the subnet $\{S_{n}:n\in E\}$ in $X-U$ is $I_{E}$-convergent$(\mathcal{M})$ to $x_{0}$, by condition $\textbf{(b)}$. So $X-U\neq (X-U)^{cl}$ and hence $U$ is not open relative to $\tau$, which is a contradiction.    \\

Conversely, suppose that a net $\{P_{n}:n\in D\}$ is $I_{D}$-convergent to a point $x_{0}$ and fails to $I_{D}$-convergent$(\mathcal{M})$ to $x_{0}$. Then by condition $\textbf{(c)}$, there is a subnet $\{T_{m}:m\in E\}$ no subnet of which is ideal convergent$(\mathcal{M})$ to $x_{0}$ relative to  any nontrivial ideal. Since $I_{E}$ is $E$-admissible ideal of $E$ so by definition  for each $r\in E$ the set $B_{r}$=$\{m\in E: m\geq r\}\in F(I_{E})$. Since $\{T_{m}:m\in E\}$ is $I_{E}$-convergent to $x_{0}$ relative to $\tau$, the point $x_{0}$ must lie in the closure of each set $A_{M}$=$\{T_{m}:m\in M\}$ for each $M\in F(I_{E})$. Consequently for each $M$ in $F(I_{E})$ there is a directed set $E_{M}$ and a net $\{S(M,n):n\in E_{M}\}$ in $M$, such that the composition $\{T\circ S(M,n):n\in E_{M}\}$ lying in $A_{M}$ is $I_{E_{m}}$-convergent$(\mathcal{M})$ to $x_{0}$. Now we apply the condition $\textbf{(d)}$. If we take $R(M,f)$=$(M,f(M))$ for each $(M,f)$ in $F(I_{E})\times (\times \{E_{M}:M\in F(I_{E})\})$ then $T\circ S\circ R$ is $I_{\mathcal{F}}$-convergent$(\mathcal{M})$ to $x_{0}$, where $\mathcal{F}$=$F(I_{E})\times (\times \{E_{M}:M\in F(I_{E})\})$ and $F(I_{E})$ is directed by set inclusion '$\subset$'. Moreover for each $m\in E$ there exists $B_{m}$ in $F(I_{E})$ such that $S\circ R(B_{m},f)$=$S(B_{m},f(B_{m}))\in B_{m}$; i.e., $S\circ R(B_{m},f)\geq m$. Therefore, $T\circ S\circ R$ is a subnet of $T$ and the result follows.
\end{proof}

\textbf{Acknowledgements.}  The second author is  thankful to University Grants Commission,India for the grant of Senior Research Fellowship during the preparation of this paper. We are thankful to the referees for their valuable suggestions which improved the quality and presentation of the paper substantially.    \\

\end{document}